\pdfoutput=1

\documentclass[a4paper,11pt, german]{amsart}%
\usepackage{amsfonts}
\usepackage{amsthm}
\usepackage{amssymb}
\usepackage{graphicx}
\usepackage{float}
\usepackage[space]{cite}\usepackage{eqparbox,array}
\usepackage{algorithm}
\usepackage[noend]{algorithmic}
\usepackage{dcolumn}
\usepackage{amsmath}
\usepackage[german]{babel}
\usepackage[utf8]{inputenc}
\usepackage{epstopdf}
\usepackage{hyperref}
\usepackage{color}

\theoremstyle{plain}

\theoremstyle{definition}

\theoremstyle{plain}

\theoremstyle{remark}

\numberwithin{equation}{section}

\begin{document}
\title{Computeralgebra - vom Vorlesungsthema zum Forschungsthema}
\author[J. B\"{o}hm]{Janko B\"{o}hm}
\address{Janko B\"{o}hm\\Fachbereich Mathematik\\
Universit\"{a}t Kaiserslautern\\
Postfach 3049\\
D-67653 Kaiserslautern, Germany}
\email{boehm@mathematik.uni-kl.de}
\urladdr{\href{http://www.mathematik.uni-kl.de/~boehm/index.htm}{http://www.mathematik.uni-kl.de/~boehm/index.htm}%
}

\thanks{The author acknowledges support of SFB-TRR 195 of the German Research Foundation (DFG)}

\renewcommand{\abstractname}{Abstract}
\begin{abstract}
In this note for the joint meeting of DMV and GDM we illustrate with examples the role of computer algebra in university mathematics education. We
discuss its potential in teaching algebra, but also computer algebra as a subject in its own right, its value in the context of practical programming projects and its role as a research topic in student papers.
\end{abstract}
\maketitle

\section{Einleitung}

In diesem einführenden Text möchte ich auf die Rolle der Computeralgebra in der universitären Mathematikausbildung eingehen, insbesondere auf den Nutzen der Computeralgebra in der Grundlehre, auf die Computeralgebra als eigenständiges Lehrgebiet und Forschungsgebiet im Rahmen von Abschlussarbeiten, und die Wechselwirkung von Forschung und Lehre mit der Entwicklung von Computeralgebrasystemen. Insbesondere die Konzepte zu den Grundlagenvorlesungen in Algebra und Zahlentheorie (wie in Abschnitt \ref{sec one} beschrieben) sind auch zu erheblichen Teilen in der gymnasialen Oberstufe verwendbar, etwa in Facharbeiten, Seminaren und Leistungskursen.

\section{Grundlagenvorlesungen zur Algebra und Zahlentheorie}\label{sec one}

\begin{figure}[H]
\begin{center}
\includegraphics[height=4.75cm]{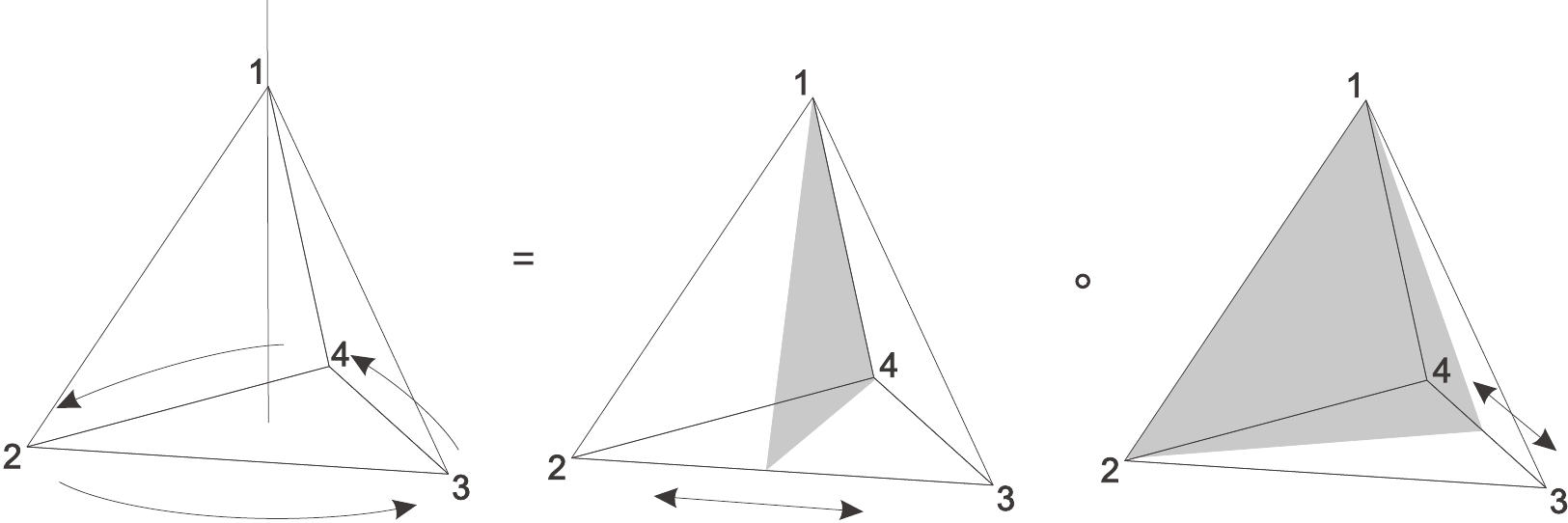}
\end{center}
 \caption{Komposition von zwei Spiegelsymmetrien des Tetraeders.}
 \label{tetraedersym}
 \end{figure}

\begin{figure}[h]
\begin{center}
\includegraphics[height=8cm]{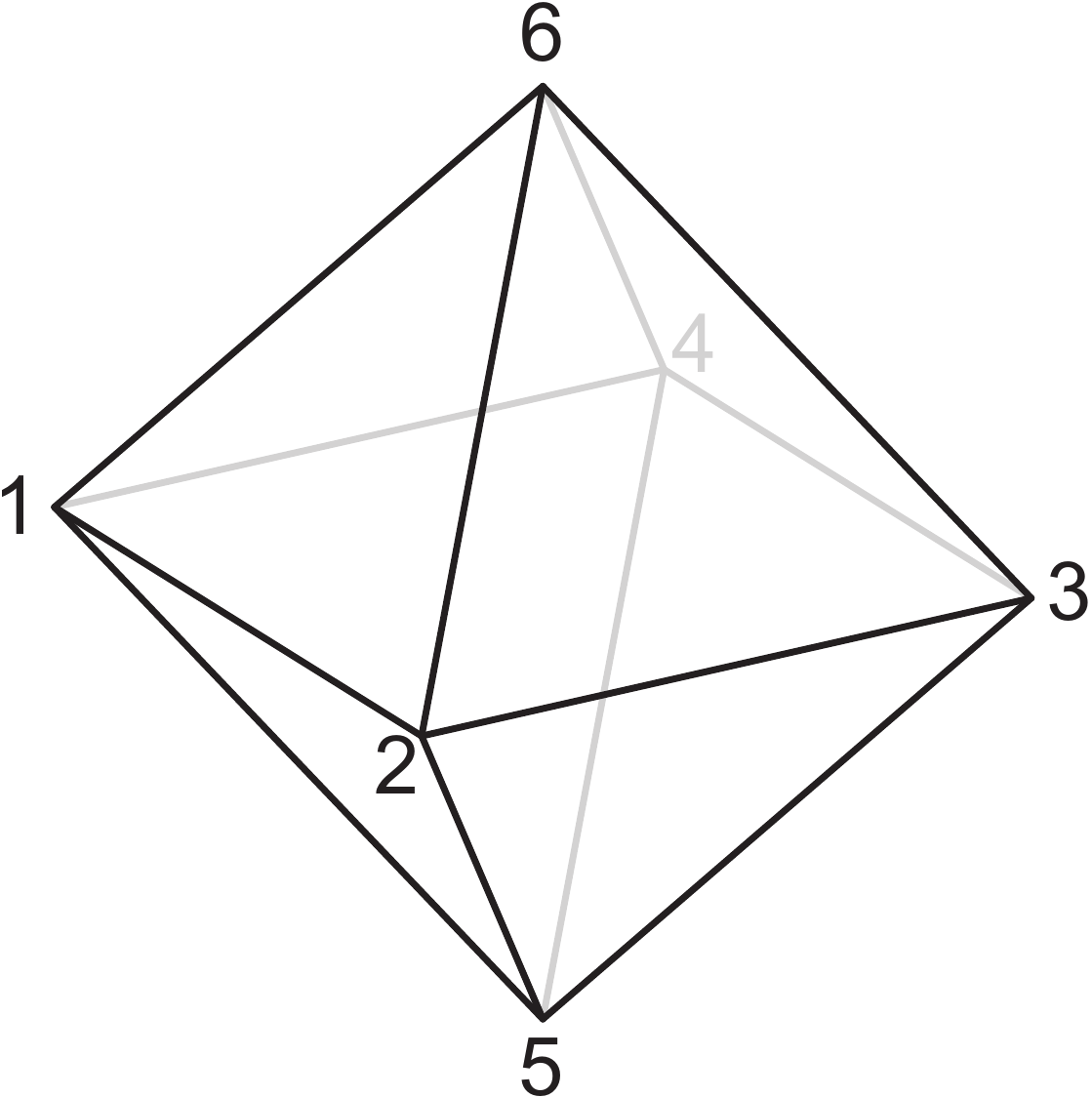}
\end{center}
 \caption{Oktaeder mit Eckennummerierung.}
 \label{oktaeder}
 \end{figure}

Im Rahmen der Grundausbildung in der Algebra und Zahlentheorie kann die Computeralgebra in einer Vielzahl von Themengebieten sowohl in Vorlesungen als auch in Übungen verwendet werden, um abstrakte Konzepte explizit und anschaulich zu erklären und Experimente zu ermöglichen. Beispiele sind etwa der Euklidische Algorithmus, die Lösung von simultanen Kongruenzen, die Illustration des Gruppenbegriffs inklusive weiterführender Konzepte wie der Sylowsätze, die lineare Algebra über Euklidischen Ringen mit Anwendung auf die Jordansche Normalform, und Restklassenringe und ihre Rolle bei Primzahltests und Public-Key-Kryptographie. Siehe dazu etwa die Anwendungsbeispiele in \cite{book}. Zum Beispiel lässt sich das Gruppenkonzept anhand von Symmetriegruppen von Platonischen Körpern einführen. Wie in Abbildung~\ref{tetraedersym} für die Hintereinanderausführung von zwei Spiegelungen des Tetraeders illustriert, ist die Komposition von zwei Symmetrien wieder eine Symmetrie.
Weiter gibt es zu jeder Symmetrie eine inverse Symmetrie, die Symmetrien bilden also eine Gruppe, und deren Elemente können als Permutationen der Ecken dargestellt werden. Die Symmetriegruppe des Tetraeders lässt sich leicht als die symmetrische Gruppe $S_4$ identifizieren, d.h. jede Permutation der Ecken ist erlaubt. Bei dem Oktaeder in Abbildung \ref{oktaeder} ist die Situation etwas komplizierter.
 
Während sich leicht Symmetrien des Oktaeders finden lassen, bleibt das Problem, alle Symmetrien aufzuzählen und zu zeigen, dass es keine weiteren gibt. Mit Hilfe der Bahnenformel kann man die Anzahl der Symmetrien bestimmen, in diesem Fall sind es $48$. Haben wir also einige Symmetrien gefunden, so können wir die davon erzeugte Untergruppe betrachten. Mit Hilfe des Computeralgebrasystems \textsc{GAP} \cite{gap}, bestimmen wir die Ordnung dieser Untergruppe, z.B. für die Spiegelung am Äquator, eine $90$-Grad-Drehung um die dazu senkrechte Achse und eine weitere solche Drehung:

\noindent\texttt{\color{blue}gap> \color{black}G:=Group((5,6), (1,2,3,4), (2,5,4,6));}

\noindent\texttt{\color{blue}gap> \color{black}Size(G);}

\noindent\color{red}\texttt{48}\color{black}

\noindent Wir können also folgern, dass $G$ schon die Symmetriegruppe des Oktaeders ist und nun z.B. leicht mit Hilfe von \textsc{GAP} alle Elemente von $G$ aufzählen und $G$ genauer analysieren. Weitergehende Konzepte, die sich auf ähnliche Weise vermitteln lassen, sind etwa die Quotientengruppe und der Homomorphiesatz. Für mehr Details siehe z.B. \cite[Abschnitt 2.2]{book}.

\section{Der erste Kontakt mit der Computeralgebra als Vorlesungsthema}\label{sec two}

\begin{figure}[h]
\begin{center}
\includegraphics[width=14cm]{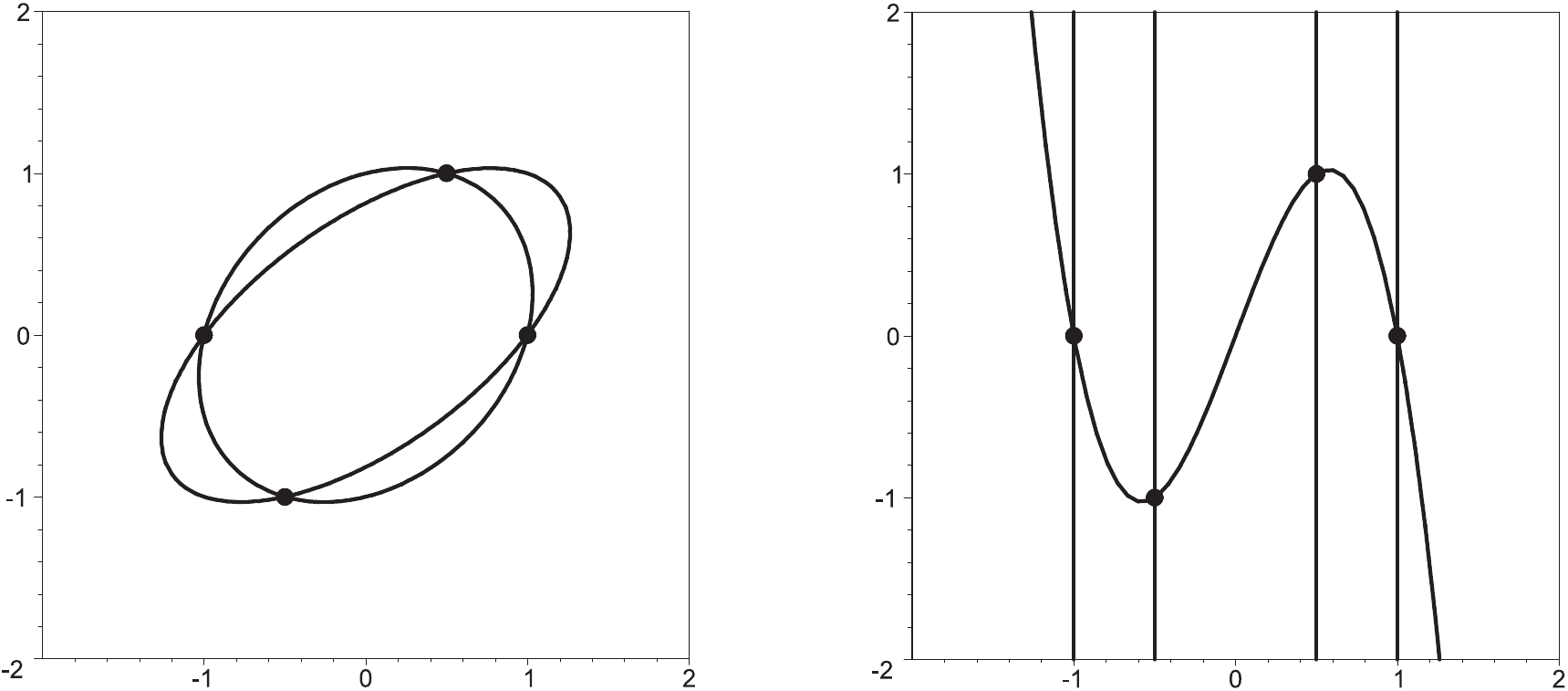}
\end{center}
 \caption{Elimination mit dem Buchbergeralgorithmus.}
 \label{buchberger}
 \end{figure}

Erstmals als eigenständiges Vorlesungsthema tritt die Computeralgebra an der TU Kaiserslautern im Bacherlorstudiengang Mathematik im Rahmen der Vorlesung \emph{Symbolisches Rechnen} auf, einer von vier Veranstaltungen zur praktischen Mathematik. Zentrale Themen sind hier Faktorisierungsalgorithmen und die Verallgemeinerung des Gaußalgorithmus zur Lösung von linearen Gleichungssystemen auf nichtlineare polynomiale Systeme. Der wesentliche Algorithmus ist in diesem Zusammenhang der Buchbergeralgorithmus zur  Bestimmung einer Gröbnerbasis. Wie beim Gaußalgorithmus können hiermit Variablen eliminiert werden, Abbildung \ref{buchberger} zeigt eine solche Transformation.

Das Computeralgebrasystem \textsc{Singular}, das an der TU Kaiserslautern entwickelt wird \cite{Singular} ist eines der führenden Systeme zur Behandlung von polynomialen Gleichungssystemen. Der Durchschnitt der beiden Ellipsen lässt sich damit wie folgt bestimmen:

\noindent\texttt{\color{blue}> \color{black}ring R=0,(y,x),lp;}

\noindent\texttt{\color{blue}> \color{black}ideal I = 2x2-xy+2y2-2, 2x2-3xy+3y2-2;}

\noindent\texttt{\color{blue}> \color{black}groebner(I);}

\noindent\texttt{\color{red}\_[1]=4x4-5x2+1}

\noindent\texttt{\color{red}\_[2]=3y+8x3-8x}

\noindent Wir beobachten, dass die beiden Ellipsengleichungen in ein äquivalentes System transformiert werden, das aus einer univariaten Gleichung und einer in der anderen Variablen linearen Gleichung besteht. Die Lösungen der univariaten Gleichung lassen sich dann leicht mit symbolischen oder numerischen Methoden bestimmen und in die andere Gleichung einsetzen. Für mehr Details siehe z.B. \cite{cox}.

\section{Programmierpraxis in der Computeralgebra}

\begin{figure}[H]
\begin{center}
\includegraphics[height=7cm]{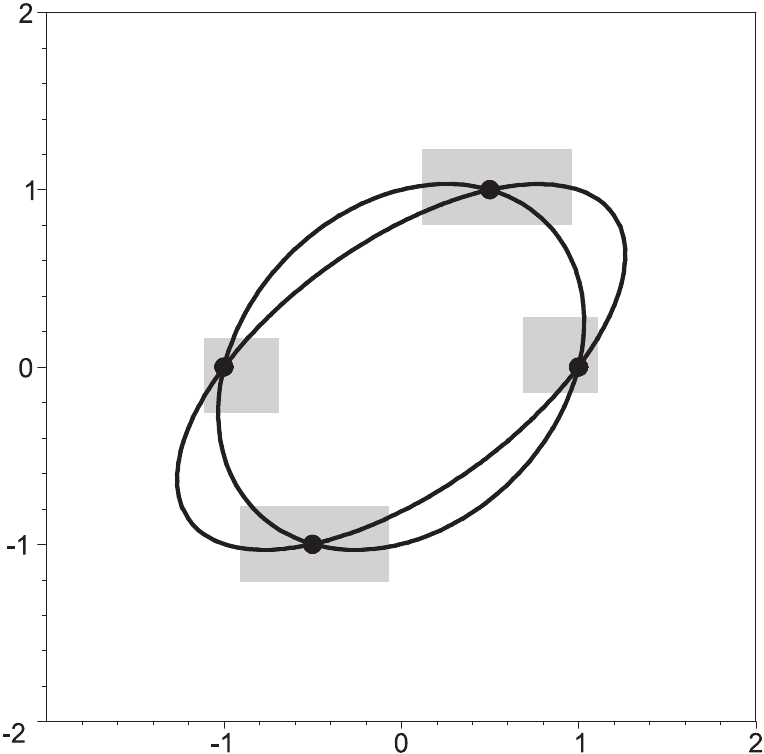}
\end{center}
 \caption{Rootisolation f\"ur den Durchschnitt von zwei Ellipsen.}
 \label{rootisolation}
 \end{figure}

Praktische Programmiererfahrung sammeln die Studierenden an der TU Kaiserslautern im Bachelorstudiengang im Rahmen eines Fachpraktikums. In dieser Projektarbeit wird ein mathematisches Problem algorithmisch behandelt. Ein Beispiel für ein solches Projekt ist etwa die Implementation eines Verfahrens zur Isolation von Lösungen eines polynomialen Gleichungssystems in Boxen mittels Intervallarithmetik, siehe dazu Abbildung \ref{rootisolation}. In jeder Box soll hier jeweils genau eine L\"osung liegen.

\section{Abschlussarbeit im Bachelorstudiengang}

\begin{figure}[h]
\begin{center}
\includegraphics[width=14cm]{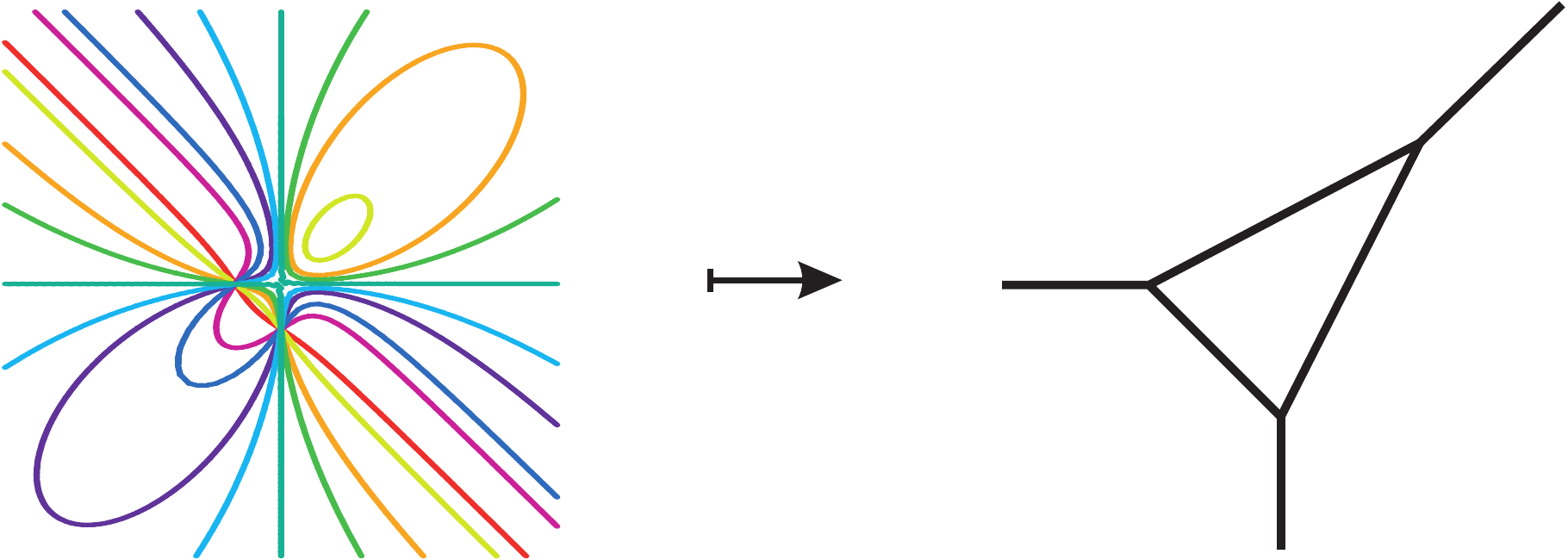}
\end{center}
 \caption{Tropikalisierung einer Familie von elliptischen Kurven.}
 \label{tropical}
 \end{figure}

Die erwähnten Praktika führen die Studierenden oft direkt zu einem Thema für eine Bachelorarbeit. Hier beginnt in manchen Fällen schon die Beschäftigung mit wissenschaftlichem Neuland. Ein Beispiel ist eine Arbeit zur massiv-parallelen Berechnung von tropischen Varietäten \cite{Bendle}. Diese Varietäten stellen eine kombinatorische Version von Lösungsmengen von polynomialen Systemen dar, siehe Abbildung \ref{tropical} f\"ur die Tropikalisierung einer Familie von algebraischen elliptischen Kurven. Die genannte Arbeit ist Teil eines umfangreichen Projekts zur Verwendung von \textsc{Singular} auf großen Clustern \cite{BDFPRR}.

\section{Weiterführende Vorlesungen in Algebra und Geometrie}

Die Computeralgebra spielt ebenso eine wichtige Rolle in den Masterstudiengängen mit Vertiefung in Algebra oder Geometrie. Die eigentliche Vorlesung \emph{Computeralgebra} behandelt in Kaiserslautern z.B. Themen wie Primärzerlegung (das geometrische Analogon zur Primfaktorisierung von Zahlen), die lineare Algebra über Polynomringen und die lokale Geometrie von algebraischen Mengen. Eng verknüpft damit sind Vorlesungen zur kommutativen Algebra, zur algebraischen Geometrie, zur algorithmischen Zahlentheorie, zur tropischen Geometrie und zur Singularitätentheorie.

\section{Forschung und Lehre als Fundament der Entwicklung von Computeralgebrasystemen}

Abschließend ist zu bemerken, dass nicht nur die Lehre von der Computeralgebra profitiert, sondern gerade die Arbeit an Algorithmen und Implementierungen in Bachelor-, Master- und Doktorarbeiten eine wesentliche Rolle in der Weiterentwicklung von Computeralgebrasystemen wie \textsc{Singular} spielt. Ein aktuelles Projekt ist die Entwicklung eines neuen Open Source Computeralgebrasystems \textsc{OSCAR} im Rahmen des Sonderforschungsbereichs TRR 195, das auf den bestehenden Systemen aufbaut und kommutative Algebra bzw. algebraische Geometrie, Gruppentheorie, polyhedrale Geometrie und Zahlentheorie zusammenbringt.

\end{document}